\documentclass[11pt]{article}

\usepackage[a4paper,margin=1in]{geometry}
\usepackage{amsmath,amssymb,amsthm}
\usepackage[hidelinks]{hyperref}
\usepackage{enumitem}
\usepackage{tikz}
\usepackage{needspace}

\usepackage{float}
\newtheorem{theorem}{Theorem}
\newtheorem{lemma}[theorem]{Lemma}
\newtheorem{corollary}[theorem]{Corollary}
\theoremstyle{definition}
\newtheorem{definition}[theorem]{Definition}
\theoremstyle{remark}

\newcommand{\Kcal}{\mathcal K}
\newcommand{\Wcal}{\mathcal W}

\title{Contractable 2-connected subgraphs in minimally 3-connected graphs}
\author{\vspace{2mm} Rongchuan Tao\thanks{Corresponding author. E-mail: rctao@hku.hk.} 
\\ Department of Mathematics\\University of Hong Kong \\ Hong Kong, China\smallskip\\}
\date{}

\begin{document}
\maketitle

\begin{abstract}
We study the 2-connected subgraphs of minimally 3-connected graphs whose contraction preserves 3-connectedness. We show that every minimally 3-connected graph that is not a wheel or a complete bipartite graph $K_{3,n}$ contains a contractable 2-connected induced subgraph that is a cycle, a fan, a $K_{2,n}$ or a $K_{3,n}$. In addition, these subgraphs satisfy certain degree restrictions. We show that none of these four families can be discarded since for each one of the four families, there exists a minimally 3-connected graph containing a contractable induced 2-connected subgraph in that family and none in the other three families. Reversing the contraction process gives a new construction for 3-connected graphs using edge addition, expansions into prescribed 2-connected subgraphs and $Y$-sum.

\vskip 4mm

\end{abstract}

\section{Introduction}\label{sec:introduction}

An edge $e=uv$ in a 3-connected graph $G$ is contractable if $G/e$ is 3-connected, or equivalently, $G-\{u,v\}$ is $2$-connected. Tutte~\cite{Tutte1961} proved that every 3-connected graph with at least five vertices has a contractable edge, and that every 3-connected graph can be constructed from wheels by edge additions and vertex splittings. Contractable edges have been studied extensively in the literature. In general, a subgraph whose contraction preserves certain property naturally corresponds to a type of chain theorem, which has applications in structural characterization of graphs. 

A cycle $C$ of a graph $G$ is non-separating if $G-V(C)$ is connected. Thomassen and Toft~\cite[Corollary 1]{ThomassenToft} proved that every $2$-connected graph of minimum degree at least three contains an induced non-separating cycle. For 3-connected graphs, they showed that either there exists an induced cycle whose contraction is 3-connected, or there exists an induced non-separating cycle with length at most 4 whose complement has at most one block containing cycles.

A connected subgraph $H$ of a 3-connected graph $G$ is contractable if $G/H$ is 3-connected. Barnette \cite[Theorem 2]{Barnette1998} proved that every 3-connected graph other than $K_5$ contains a contractable cycle, a triangle with a cubic vertex, or a $K_{2,3}$ subgraph whose three degree-two vertices are cubic in the original graph. Conversely, this reduction gives a recursive construction of 3-connected graphs in \cite[Theorem 4]{Barnette1998} from $K_4$ and $K_5$ using circuit expansions, $M$-configuration expansions and vertex splittings that produce a $t$-triangle. The contractable $M$-configuration directly corresponds to the 2-crown in our result. However, the vertex splittings alternative that produce a $t$-triangle does not require the $t$-triangle to be contractable. A natural question is: whether there is a reduction described entirely by contractable induced 2-connected subgraphs and edge deletions in a 3-connected graph?

We solve this problem affirmatively: every minimally 3-connected graph other than a wheel or a $K_{3,n}$ contains a contractable induced cycle, fan, $2$-crown, or 3-crown. All configurations are 2-connected subgraphs and the latter three configurations satisfy additional degree conditions in the original graph. Iterating the theorem, and taking a minimally 3-connected spanning subgraph after each contraction, reduces the graph to a wheel or a $K_{3,n}$. Reversing this sequence gives a construction for 3-connected graphs. 

\subsection*{Notation and main result}

All graphs are finite and simple. For any graph $G$, let $|G|:=|V(G)|$. We use the standard definitions of paths, cycles, blocks, bridges, and $k$-connected graphs~\cite{BondyMurty}. We allow a tree or a path to contain only one vertex. For $k\geq 2$, a $k$-connected graph has at least $k+1$ vertices. A vertex of degree three is called cubic. If $H$ is a nonempty connected subgraph of $G$, then $G/H$ is the graph obtained from $G$ by contracting $H$ to one vertex, deleting loops, and replacing parallel edges by single edges. A connected subgraph $H$ is \emph{contractable} if $G/H$ is 3-connected. For any vertex subset $X\subseteq V(G)$, let $N_G(X)$ be the vertices outside $X$ having a neighbour in $X$. For subgraphs $H,H'$ of $G$, let $N_{H'}(H):=N_G(V(H))\cap V(H')$ be the set of neighbours of $H$ in $H'$. Similarly define $N_{H'}(v)$ for a vertex $v$. The attachment set of a subgraph $H$ of $G$ is the neighbours of its complement in $H$, i.e. $N_{H}(G-H)$. A vertex in $G$ $meets$ (resp. $misses$) a disjoint subgraph $H$ of $G$ if it has a neighbour in $H$ (resp. has no neighbour in $H$).

\begin{definition}\label{def:objects}
Let $G$ be a graph and $H$ be an induced subgraph of $G$.
\begin{enumerate}[label=\textup{(\roman*)}]

\item $H$ is a \emph{fan} if it has vertices $x,y_1,\ldots,y_k$, where $k\ge2$, and edge set
\[
E(H)=\{xy_i:1\le i\le k\} \cup\{y_i y_{i+1}:1\le i<k\}
\]
such that $y_i$ is cubic in $G$ for any $1\leq i\leq k$. We call $x$ its centre and $y_1\cdots y_k$ its rim path.
\item $H$ is a \emph{2-crown} (resp. \emph{3-crown}) if it is isomorphic to some $K_{2,n}$ with $n\ge2$ (resp. $K_{3,n}$ with $n\ge3$), whose vertices in the colour class of size $n$ are cubic in $G$. These vertices are called the \emph{teeth} of the 2-crown or 3-crown.
\end{enumerate}
\end{definition}

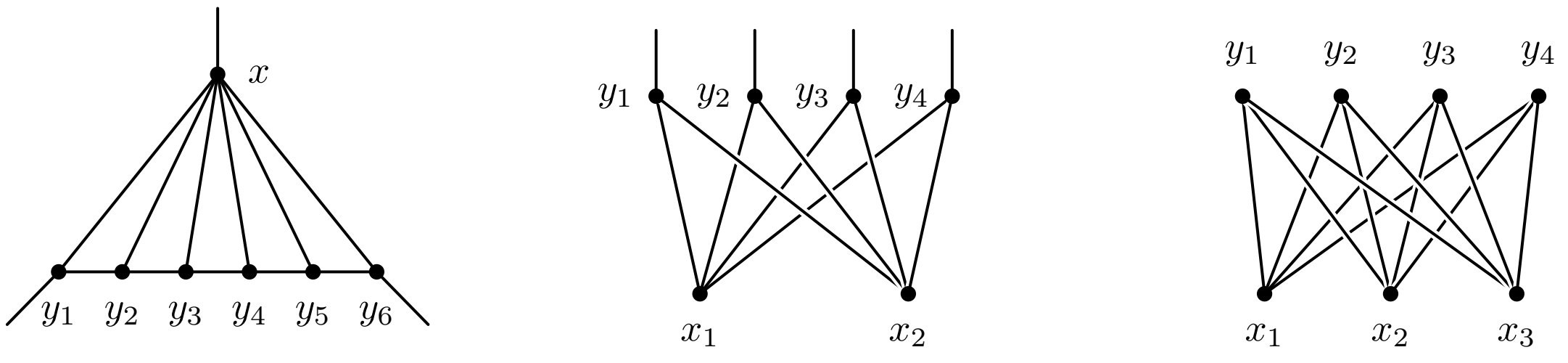
\begin{figure}[H]
  \centering
  \begingroup
  \tikzset{
    obj vertex/.style={circle,fill=black,inner sep=0pt,minimum size=4pt},
    obj edge/.style={draw=black,line width=.75pt},
    obj crossing/.style={obj edge,preaction={draw=white,line width=2.3pt}}
  }
  \begin{minipage}[t]{.32\linewidth}
  \centering
  \begin{tikzpicture}[x=1cm,y=1cm,line cap=round,line join=round,
    every node/.style={font=\small}]
    \path[use as bounding box] (-2.1,-1.65) rectangle (2.1,1.75);
    \coordinate (x) at (0,1);
    \foreach \i/\a in {1/-1.45,2/-.87,3/-.29,4/.29,5/.87,6/1.45}{
      \coordinate (y\i) at (\a,-.8);
      \draw[obj edge] (x)--(y\i);
    }
    \draw[obj edge] (y1)--(y2)--(y3)--(y4)--(y5)--(y6);
    \draw[obj edge] (x)--(0,1.6)
      (y1)--(-1.92,-1.28) (y6)--(1.92,-1.28);
    \node[obj vertex] at (x) {};
    \node[right=4pt] at (x) {$x$};
    \foreach \i in {1,...,6}{
      \node[obj vertex] at (y\i) {};
      \node[below=4pt] at (y\i) {$y_{\i}$};
    }
  \end{tikzpicture}
  \par\smallskip
  {\small (a) A fan.}
  \end{minipage}\hfill
  \begin{minipage}[t]{.32\linewidth}
  \centering
  \begin{tikzpicture}[x=1cm,y=1cm,line cap=round,line join=round,
    every node/.style={font=\small}]
    \path[use as bounding box] (-2.1,-1.65) rectangle (2.1,1.75);
    \coordinate (x1) at (-.95,-1);
    \coordinate (x2) at (.95,-1);
    \foreach \i/\a in {1/-1.35,2/-.45,3/.45,4/1.35}{
      \coordinate (y\i) at (\a,.8);
      \draw[obj edge] (y\i)--(\a,1.4);
      \draw[obj edge] (x1)--(y\i);
    }
    \foreach \i in {1,...,4}{\draw[obj crossing] (x2)--(y\i);}
    \foreach \i in {1,2}{
      \node[obj vertex] at (x\i) {};
      \node[below=4pt] at (x\i) {$x_{\i}$};
    }
    \foreach \i in {1,...,4}{
      \node[obj vertex] at (y\i) {};
      \node[left=2pt] at (y\i) {$y_{\i}$};
    }
  \end{tikzpicture}
  \par\smallskip
  {\small (b) A $2$-crown.}
  \end{minipage}\hfill
  \begin{minipage}[t]{.32\linewidth}
  \centering
  \begin{tikzpicture}[x=1cm,y=1cm,line cap=round,line join=round,
    every node/.style={font=\small}]
    \path[use as bounding box] (-2.1,-1.65) rectangle (2.1,1.75);
    \coordinate (x1) at (-1.15,-1);
    \coordinate (x2) at (0,-1);
    \coordinate (x3) at (1.15,-1);
    \foreach \i/\a in {1/-1.35,2/-.45,3/.45,4/1.35}{
      \coordinate (y\i) at (\a,.8);
      \draw[obj edge] (x1)--(y\i);
    }
    \foreach \j in {2,3}{
      \foreach \i in {1,...,4}{\draw[obj crossing] (x\j)--(y\i);}
    }
    \foreach \i in {1,...,3}{
      \node[obj vertex] at (x\i) {};
      \node[below=4pt] at (x\i) {$x_{\i}$};
    }
    \foreach \i in {1,...,4}{
      \node[obj vertex] at (y\i) {};
      \node[above=4pt] at (y\i) {$y_{\i}$};
    }
  \end{tikzpicture}
  \par\smallskip
  {\small (c) A 3-crown.}
  \end{minipage}
  \endgroup
  \caption{Examples of the configurations in Definition~\ref{def:objects}. There may exist additional edges incident with $x$ or some $x_i$.}
  \label{fig:defined-objects}
  \end{figure}

See Figure 1 for examples of fan, 2-crown and 3-crown. A fan with $k=2$ is a triangle and hence an induced cycle; we include this case for convenience. Let $\Kcal$ be the class of graphs obtained from $K_{3,n}$, $n\ge3$, by adding any set of edges among the three vertices in its part of size three, and let $\Wcal$ be the class of wheels, including $K_4$. A diamond is the graph obtained by adding an edge to a pair of nonadjacent vertices in a cycle of length four. The following theorem is our main result.

\begin{theorem}\label{thm:main}
Let $G$ be a minimally 3-connected graph. If $G$ is neither a wheel nor a graph $K_{3,m}$ with $m\ge3$, then $G$ has a contractable induced subgraph that is either a cycle, a fan, a 2-crown or a 3-crown.
\end{theorem}

Note that wheels and the graphs $K_{3,m}$ have no contractable $2$-connected subgraph, so they are irreducible under contraction of $2$-connected subgraphs. The four alternatives in Theorem~\ref{thm:main} are individually necessary. For each alternative, each of the following examples gives a minimally 3-connected graph with a contractable induced subgraph of that type and no contractable induced subgraph of any of the other three types. 

\begin{enumerate}[label=\textup{(\arabic*)},leftmargin=*]
\item \emph{Cycles.}
The Petersen graph is cubic and minimally
3-connected. Its contractable induced $2$-connected subgraphs are precisely its $5$-cycles.

\item \emph{Fans.}
Let $D,D'$ be disjoint diamonds with $V(D)=\{a,b,c,d\}$ and $V(D')=\{a',b',c',d'\}$ such that $abc$, $abd$, $a'b'c'$ and $a'b'd'$ are triangles. Add the edges $bb',cc',dd'$. This graph is minimally 3-connected. Among the four types, its contractable induced subgraphs are precisely the two diamonds, and each is a fan.

\item \emph{$2$-crowns.}
Take two disjoint copies of $K_{2,3}$ with bipartition $(X,A)$ and $(Y,B)$, where $X=\{x_1,x_2\}$, $Y=\{y_1,y_2\}$, $A=\{a_1,a_2,a_3\}$ and $B=\{b_1,b_2,b_3\}$. Add edges $a_iy_1$, $b_ix_1$ for each $1\leq i\leq 3$
and $x_2y_2$. The graph is minimally 3-connected. Among the four types, its
contractable induced subgraphs are only four $2$-crowns.

\item \emph{3-crowns.}
Take two disjoint copies of $K_{3,4}$ with bipartitions $(X,A)$ and
$(Y,B)$, where $X=\{x_1,x_2,x_3\}$ and
$Y=\{y_1,y_2,y_3\}$, and add edge $x_i y_i$ for $1\le i\le3$. This graph is minimally 3-connected and among the four types its contractable induced subgraphs are precisely the two $K_{3,4}$ subgraphs, both 3-crowns.

\end{enumerate}

The outline of the proof is as follows. We start from a minimally 3-connected graph $G$ that is not a wheel or the $K_{3,n}$. Dirac's theorem~\cite{Dirac1963} says the minimally 3-connected graphs with no two vertex-disjoint cycles are exactly wheels and the $K_{3,n}$'s. Then by the non-separating induced-cycle theorem of Thomassen and Toft~\cite{ThomassenToft}, either the graph $G$ has a contractable induced cycle, or it contains a non-separating cycle $C$ of length at most four such that $G-C$ has exactly one block that contains a cycle. We analyze the trees attached to that block and show that certain paths in the tree give contractable fans, while certain stars in the tree give contractable 2-crowns or 3-crowns. Finally, if none of these paths or stars appear, then every component of $G-C-B$ must be a path, and a comparison of their attachments produces a contractable crown or fan.

\section{Preliminaries}\label{sec:connectivity}

We first recall Halin's lemma: every cycle of a minimally 3-connected graph contains at least two cubic vertices \cite[Satz~5, p.~184]{Halin1969Minimal}. The following consequence recognizes induced fans and $2$-crowns from prescribed edges.

\begin{lemma}\label{lem:forms}
Let $G$ be a minimally 3-connected graph with $|G|\ge5$.
\begin{enumerate}[label=\textup{(\arabic*)},leftmargin=*]
\item Suppose $S=\{x,y_1,\ldots,y_k\}\subsetneq V(G)$ with $k\ge2$, and $xy_i$, $y_i y_{i+1}$ are all edges. Then $S$ induces a fan in $G$, whose centre need not be $x$.
\item Suppose $S=\{x_1,x_2,y_1,\ldots,y_k\}\subseteq V(G)$ with $k\ge2$, each $y_j$ is cubic and each $x_iy_j$ is an edge. Then $S$ induces a $2$-crown or a fan.
\end{enumerate}
\end{lemma}
\begin{proof}
\begin{enumerate}[label=\textup{(\arabic*)},leftmargin=*]
\item If $k=2$, the induced subgraph is a triangle. Thus $S$ induces a fan. 

If $k=3$, then apply Halin's lemma on $xy_1y_2$ and $xy_2y_3$, and notice that $G$ is 3-connected. We conclude that either $y_1,y_2,y_3$ are cubic or $x,y_1,y_3$ are cubic. Since $G$ is 3-connected, $y_1y_k$ is not an edge. Thus $S$ induces a fan whose centre is either $x$ or $y_2$. 
  
If $k\ge4$, then $d_G(x)\ge4$. Apply Halin's lemma on each $xy_i y_{i+1}$, we have $y_i$ is cubic for all $i$. Since $G$ is 3-connected, $y_1y_k$ is not an edge. It follows that $S$ induces a fan with centre $x$, proving the assertion.

\item The set $Y=\{y_1,...,y_k\}$ is independent since otherwise $\{x_1,x_2\}$ is a 2-vertex cut. If $x_1x_2\notin E(G)$, then $S$ induces a 2-crown. If $x_1x_2\in E(G)$, then $|Y|\leq 2$ since otherwise $x_1$, $x_2$ have degree at least 4, contrary to Halin's lemma on the triangle $x_1x_2y_1$. Now $|Y|\leq 2$ implies $S$ induces a diamond. The assertion follows from (1).
\end{enumerate}
\end{proof}

\begin{lemma}\label{lem:contract}
If $G$ is 3-connected, $H$ is connected, and
$|V(G-H)|\ge3$, then $H$ is contractable if and only if
$G-V(H)$ is $2$-connected.
\end{lemma}
\begin{proof}
If $H$ is contractable, then $G/H$ is 3-connected and hence $G-V(H)$ is $2$-connected. Conversely, suppose $G-V(H)$ is $2$-connected, then deleting two vertices of $G/H$ preserves connectivity: if the contracted vertex is deleted, then the $2$-connectivity of $G-V(H)$ implies the connectivity. If this is not the case, then the 3-connectivity of $G$ implies the connectivity. It follows that $G/H$ is 3-connected and $H$ is contractable.
\end{proof}

\begin{lemma}\label{lem:treepath}
Let $R_1,R_2$ be vertex-disjoint trees with $|R_1|+|R_2|\ge3$. Regard the unique vertex of a singleton tree as a leaf. If every leaf of $R_1$ has a neighbour on $R_2$ and every leaf of $R_2$ has a neighbour in $R_1$, then their union with these joining edges is $2$-connected.
\end{lemma}
\begin{proof}
The statement is obvious if one of $R_1$ and $R_2$ is a singleton. Thus we may assume both trees are nontrivial. Suppose $v\in V(R_1)$, then each component of $R_1-v$ contains a leaf of $R_1$ and therefore joins $R_2$. A similar argument applies to $R_2$. It follows that the resulting graph is 2-connected. 
\end{proof}

We will apply this lemma with one of the trees a path. The following lemma will be used to deal with the unique block $B$ in the complement of a nonseparating cycle.

\Needspace{5\baselineskip}
\begin{lemma}\label{lem:block}
The following statements hold.
\begin{enumerate}[label=\textup{(\roman*)},leftmargin=*]
\item If $B$ is a 2-connected subgraph of $J$ such that $J/B$ is 2-connected, then $J$ is 2-connected if and only if at least two vertices of $B$ meet $J-V(B)$.

\item Suppose $G$ is 3-connected, $B,H$ are disjoint 2-connected subgraphs, and $w\in V(B)$. If every edge from $B-w$ to $G-V(B)$ ends in $H$, then $G[(V(B)-\{w\})\cup V(H)]$ is $2$-connected.
\end{enumerate}
\end{lemma}
\begin{proof}
For (i), if at most one vertex of $B$ meets $J-V(B)$, then $J$ has a cut vertex. If at least two vertices of $B$ meet $J-V(B)$, and suppose $v$ is a vertex cut of $J$. Then $v\in V(B)$ since $J/B$ is 2-connected. Since $B-v$ is connected and has an attachment to the connected graph $J-V(B)$, $J-v$ is connected and a contradiction.

For (ii), let $J=G[(V(B)-\{w\})\cup V(H)]$.
Then $B-w$ is connected and at least two vertices of $H$ have neighbours in it: otherwise deleting $w$ and at most one vertex separates $B-w$ from the rest of $H$. Thus $J-v$ is connected for any $v\in V(H)$. For $v\in V(B)-\{w\}$, any component of $J-v$ intersects $H$ since otherwise $\{w,v\}$ is a 2-vertex cut of $G$. Thus $J$ is 2-connected.
\end{proof}

The following two results fix the starting configuration. The first is Dirac's theorem about two disjoint cycles in a 3-connected graph.

\begin{theorem}[Dirac {\cite[Theorem~1]{Dirac1963}}]\label{thm:dirac}
A 3-connected graph has no two vertex-disjoint cycles if and only if it is a wheel, $K_5$, $K_5-e$, or a member of $\Kcal$.
\end{theorem}

In particular, a minimally 3-connected graph other than a wheel or a graph $K_{3,m}$ contains two vertex-disjoint cycles.

\begin{theorem}[Thomassen--Toft
{\cite[Theorem~2]{ThomassenToft}}]\label{lem:core}
Let $G$ be 3-connected. There is an induced cycle $C$ such that $F:=G-V(C)$ is connected and has a block $B$ containing every cycle of $F$. Moreover, one of the following holds.
\begin{enumerate}[label=\textup{(\alph*)},leftmargin=*]
\item $F=B$.
\item An induced $K_{2,3}$ contains $C$; deleting this $K_{2,3}$
leaves a connected graph containing $B$; and its three degree-two vertices have degree three in $G$.
\item An induced diamond contains $C$; deleting this diamond leaves a connected graph containing $B$; and its two degree-two
vertices have degree three in $G$.
\item A $K_4$ contains $C$, its deletion is $B$, and all but possibly one of its vertices have degree four in $G$; the possible
exception has degree three.
\end{enumerate}
\end{theorem}

We use one feature of the choice in their proof \cite[pp.~206--207]{ThomassenToft}. They first maximize the order of a block in the complement of an induced cycle. Keeping this block $B$ fixed, they then maximize the component containing $B$ in the complement of an induced cycle disjoint from $B$; exchanges of the cycle give the four alternatives above. If $G$ contains two vertex-disjoint cycles, then the block whose order is maximized contains a cycle, which means $B$ can be chosen $2$-connected. We use the theorem with this choice throughout.

For every subgraph $H$ with $3\le |V(H)|<|V(G)|$, we have
\begin{equation}\label{eq:attachments}
\bigl|\{v\in V(H):N_G(v)\setminus V(H)\ne\varnothing\}\bigr|\ge3
\end{equation}
since otherwise deleting at most two vertices of $H$ would separate the remaining vertices of $H$ from the rest of the graph. 

For the rest of this paper, let $G$ be a minimally 3-connected graph that is neither a wheel nor a $K_{3,n}$. Then $G$ has two vertex-disjoint cycles and we can choose $C,B$ as above. If \textup{(a)} is the case, then $F=B$ and $C$ is a contractable induced cycle. Alternative~\textup{(d)} is excluded by Halin's lemma, since every triangle in that $K_4$ has at most one cubic vertex. Thus we may assume either \textup{(b)} or~\textup{(c)} is the case.

Every block of $F$ other than $B$ is an edge, so contracting $B$ to a vertex $b$ produces a tree $T$. Note that $B$ is contractable in $G$ since $G-B$ is 2-connected: any vertex of $G-B$ has two internally disjoint paths to $C$. Vertices of $T-b$ are identified with the corresponding vertices of $F-B$. 

\begin{lemma}\label{lem:leaves}
If either \textup{(b)} or~\textup{(c)} is the case, then every leaf of $T$ distinct
from $b$ is cubic and has exactly two neighbours on $C$.
\end{lemma}
\begin{proof}
A leaf $v\ne b$ has one neighbour in $F$ and therefore at least two on $C$. Suppose it has at least three neighbours on $C$. Since either \textup{(b)} or~\textup{(c)} is the case, the cycle $C$ has length at most 4. For every edge $xy$ of $C$ whose ends are adjacent to $v$, both $x,y$ are cubic by Halin's lemma. Therefore, if $|C|=3$, then all its vertices are cubic, and deleting $v$ separates $C$ from $F-v$. If $|C|=4$, then at least three of its vertices are cubic; deleting $v$ and the possible fourth vertex separates those three vertices from $F-v$. Both contradict 3-connectivity.
\end{proof}

We also record the corresponding observation at $b$:
\begin{equation}\label{eq:root-neighbours}
 d_T(b)+|N_C(B)|\ge3.
\end{equation}
Indeed, $B$ has $d_T(b)$ neighbours in $F-B$. If the sum were at most two, these neighbours together with $N_C(B)$ would separate $B$ from the remaining vertices of $C$. In particular, if $b$ is a leaf of $T$, it has at least two neighbours on $C$ in $G/B$; if $d_T(b)=2$, it has at least one.

\begin{lemma}\label{lem:triangle}
If ~\textup{(c)} is the case, then $G$ has a contractable fan.
\end{lemma}
\begin{proof}
Write $C=pqrp$, and let $u$ be the fourth vertex of the diamond, adjacent to $p,q$. The vertices $r,u$ are cubic, and $u\ne b$ is a leaf of $T$. Halin's lemma on $pqup$ allows us to assume that $q$ is cubic, so $N_G(q)=\{p,r,u\}$. Let $t$ be the unique $F$-neighbour of $r$.

By Lemma~\ref{lem:leaves}, every leaf of $T$ other than $b,u$ has cycle neighbours $p,r$, so it must be $t$. If $b$ is a leaf, then both $p,r$ meet $B$ by \eqref{eq:root-neighbours}, hence $t\in B$. Thus $T$ has exactly two leaves and is a path, and adjoining $q,r$ to it gives a cycle in $G/B$.

Deleting $b$ from this cycle gives a path $P$ in $G$. Every vertex $v\ne b$ with $d_T(v)=2$ has two tree neighbours and can meet only $p$ on $C$, so it is cubic and adjacent to $p$. The same conclusion holds for the remaining vertices of $T-b$ by Lemma~\ref{lem:leaves}; the vertices $q,r$ are also cubic and adjacent to $p$. Consequently $G[V(P)\cup\{p\}]$ is an induced fan with cubic rim $P$. Its complement is $B$, so it is contractable by Lemma~\ref{lem:contract}.
\end{proof}

\section{Proof of the main theorem}

It remains to consider the case (b). Then $C$ has length 4 and we suppose $C=prqsp$, where a leaf $u\ne b$ of $T$ is adjacent to $p,q$ with $r,s,u$ cubic. 

For any vertex $v\neq b$ of $T$, let $Z_v\subseteq\{r,s\}$ consist of the vertices that have unique $F$-neighbour $v$. Suppose $e\in E(T)$ and $R$ is a component of $T-e$. A leaf of $R$ is \emph{$R$-old} if it is a leaf of $T$ and otherwise \emph{$R$-new}. Note that there is at most one $R$-new leaf. If $v$ is an $R$-new leaf, then $v$ is the end of $e$ that belongs to $R$, and $d_T(v)=2$.

Every $R$-old leaf has at least two neighbours on $C$ in $G/B$, and every $R$-new leaf has at least one. For vertices different from $b$, this follows from Lemma~\ref{lem:leaves} and minimum degree; for $b$, it follows from \eqref{eq:root-neighbours}.

For each leaf $v\ne b$ of $T$, let $P_v=v_0v_1\cdots v_m$ with $v_0=v$ be the maximal path starting at $v$ such that $v_i\ne b$ and $d_T(v_i)=2$ for $1\le i\le m$. Thus the path stops before $b$ or a vertex of degree at least three in $T$, and it may consist only of $v$. We construct a \emph{terminal fan} in either of the following
two ways (see Figure 2 for examples of the construction).
\begin{enumerate}[label=\textup{(\roman*)},leftmargin=*]
\item If $N_C(v)=\{x,y\}$, where $x\in\{p,q\}$ and $y\in\{r,s\}$, let $j$ be largest such that every vertex of $S=v_0\cdots v_j$ is adjacent to $x$. Take $H=G[V(S)\cup\{x,y\}]$, with centre $x$.
\item If $N_C(v)=\{p,q\}$, $m\ge1$, and $v_1$ is adjacent to $x\in\{p,q\}$, choose $S=v_0\cdots v_j$ in the same way and take $H=G[V(S)\cup\{x\}]$, with centre $x$. Note that $j\geq 1$ in this case.
\end{enumerate}

\begingroup
\tikzset{
tf edge/.style={draw=black,line width=.75pt},
tf selected edge/.style={draw=black,line width=1.05pt},
tf vertex/.style={circle,fill=black,inner sep=0pt,minimum size=4pt},
tf selected/.style={tf vertex},
tf continuation/.style={draw=black,line width=.75pt,
  dash pattern=on 2.5pt off 2pt},
tf remainder/.style={draw=black,line width=.55pt,
  dash pattern=on 3pt off 2.2pt},
tf label/.style={font=\small,inner sep=2pt},
tf title/.style={font=\small,inner sep=0pt}
}
\begin{figure}[H]
\centering
\begin{minipage}[t]{.48\linewidth}
\centering
\begin{tikzpicture}[x=1cm,y=1cm,font=\small]
  \path[use as bounding box] (-2.65,-2.95) rectangle (2.65,2.65);
  \coordinate (p) at (-2,0);
  \coordinate (r) at (0,2);
  \coordinate (q) at (2,0);
  \coordinate (s) at (0,-2);
  \coordinate (v0) at (0,-1.25);
  \coordinate (v1) at (0,-.55);
  \coordinate (vj) at (0,.15);

  \path[tf remainder] (0,1.08) circle[radius=.60];
  \node[tf label,text=black] at (0,1.08) {$R$};
  \draw[tf edge] (0,0) circle[radius=2];
  \draw[tf selected edge] (p) arc[start angle=180,end angle=270,radius=2];
  \draw[tf selected edge] (s)--(v0)--(v1);
  \draw[tf selected edge,dash pattern=on 0pt off 3pt,line cap=round] (v1)--(vj);
  \draw[tf selected edge] (p)--(v0) (p)--(v1) (p)--(vj);
  \draw[tf continuation] (vj)--(0,.48);

  \node[tf selected] at (p) {};
  \node[tf vertex] at (r) {};
  \node[tf vertex] at (q) {};
  \node[tf selected] at (s) {};
  \foreach \v in {v0,v1,vj}{\node[tf selected] at (\v) {};}
  \node[left=3pt] at (p) {$p$};
  \node[above=3pt] at (r) {$r$};
  \node[right=3pt] at (q) {$q$};
  \node[below=3pt] at (s) {$s$};
  \node[tf label,anchor=north west] at (.06,-1.27) {$v_0$};
  \node[tf label,anchor=west] at (.08,-.55) {$v_1$};
  \node[tf label,anchor=west] at (.08,.15) {$v_j$};
  \node[tf label,text=black] at (-.26,-.77) {$S$};
  \node at (-1.7,1.8) {$C$};
  \node[tf title] at (0,-2.82) {\textup{(i)}\quad $N_C(v_0)=\{p,s\}$};
\end{tikzpicture}
\end{minipage}\hfill
\begin{minipage}[t]{.48\linewidth}
\centering
\begin{tikzpicture}[x=1cm,y=1cm,font=\small]
  \path[use as bounding box] (-2.65,-2.95) rectangle (2.65,2.65);
  \coordinate (p) at (-2,0);
  \coordinate (r) at (0,2);
  \coordinate (q) at (2,0);
  \coordinate (s) at (0,-2);
  \coordinate (v0) at (0,-1.25);
  \coordinate (v1) at (0,-.55);
  \coordinate (vj) at (0,.15);

  \path[tf remainder] (0,1.08) circle[radius=.60];
  \node[tf label,text=black] at (0,1.08) {$R$};

  \draw[tf edge] (0,0) circle[radius=2];
  \draw[tf edge] (q)--(v0);
  \draw[tf selected edge] (v0)--(v1);
  \draw[tf selected edge,dash pattern=on 0pt off 3pt,line cap=round] (v1)--(vj);
  \draw[tf selected edge] (p)--(v0) (p)--(v1) (p)--(vj);
  \draw[tf continuation] (vj)--(0,.48);

  \node[tf selected] at (p) {};
  \foreach \v in {r,q,s}{\node[tf vertex] at (\v) {};}
  \foreach \v in {v0,v1,vj}{\node[tf selected] at (\v) {};}
  \node[left=3pt] at (p) {$p$};
  \node[above=3pt] at (r) {$r$};
  \node[right=3pt] at (q) {$q$};
  \node[below=3pt] at (s) {$s$};
  \node[tf label,anchor=north west] at (.06,-1.27) {$v_0$};
  \node[tf label,anchor=west] at (.08,-.55) {$v_1$};
  \node[tf label,anchor=west] at (.08,.15) {$v_j$};
  \node[tf label,text=black] at (-.26,-.77) {$S$};
  \node at (-1.7,1.8) {$C$};
  \node[tf title] at (0,-2.82) {\textup{(ii)}\quad $N_C(v_0)=\{p,q\}$};
\end{tikzpicture}
\end{minipage}
\caption{Two types of terminal-fan constructions with centre $p$, where $R=T-V(S)$.}
\label{fig:terminal-fans}
\end{figure}
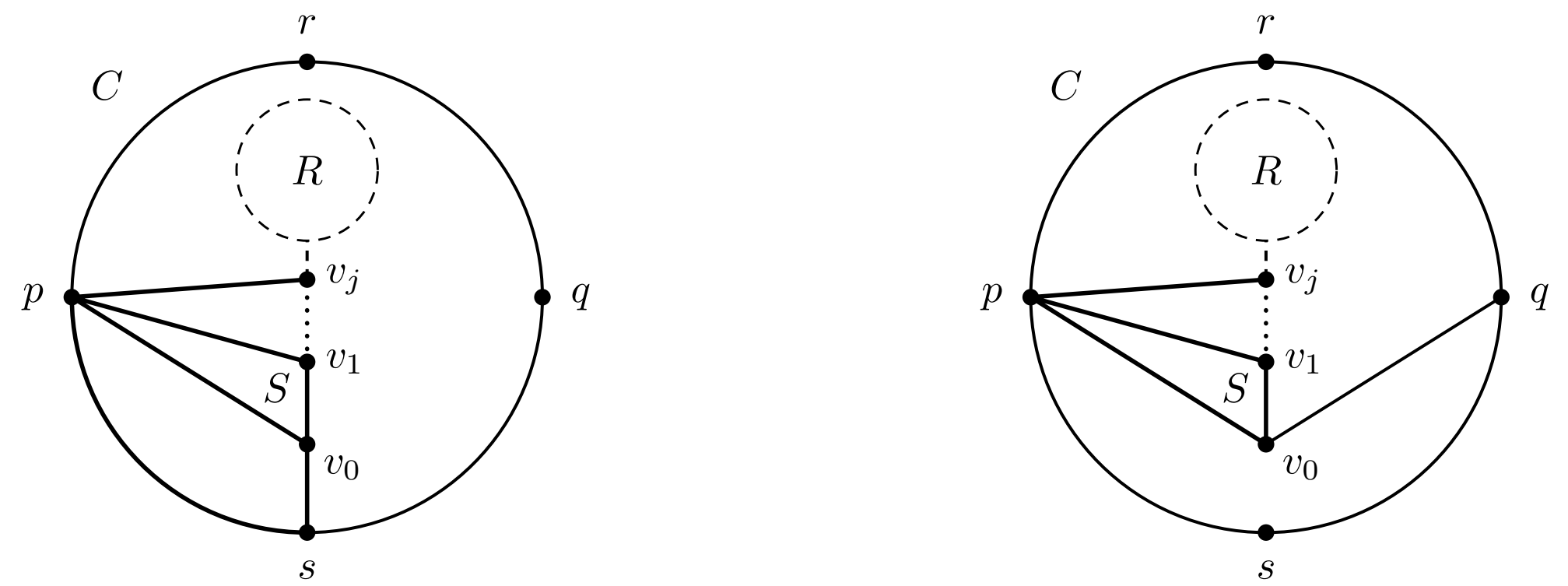
\endgroup

A terminal fan is indeed a fan with the prescribed centre by Lemma \ref{lem:forms}. A leaf defines at most one terminal fan in type~\textup{(ii)} since $v_1$ is cubic and cannot be adjacent to both $p$ and $q$. For any terminal fan $H$, define $D_H:=N_B(G-V(B)-V(H))$, the set of attachment vertices of $B$ in $G-H$.

\begin{lemma}\label{lem:fan} Let $H$ be a terminal fan defined by $S=v_0\cdots v_j$ as above. Then $D_H\ne\varnothing$ and $H$ is contractable if and only if $|D_H|\ge2$. If $H$ is not contractable, then $v_j$ has a neighbour in $B$. 
\end{lemma}

\begin{figure}[H]
\begingroup
\newcommand{\lfPanel}[1]{%
\begin{tikzpicture}[
  x=1cm,y=1cm,
  every node/.style={font=\small},
  lf edge/.style={draw=black,line width=.7pt},
  lf fan/.style={draw=black,line width=1.25pt},
  lf vertex/.style={circle,fill=black,inner sep=0pt,minimum size=4pt},
  lf Hvertex/.style={circle,fill=black,inner sep=0pt,minimum size=4pt},
  lf label/.style={inner sep=1.3pt,fill=white},
  baseline=(current bounding box.center)
]
  \path[use as bounding box] (-2.65,-3.05) rectangle (2.65,2.55);
  \coordinate (p) at (-2,0);
  \coordinate (r) at (0,2);
  \coordinate (q) at (2,0);
  \coordinate (s) at (0,-2);
  \coordinate (v0) at (0,-1.7);
  \coordinate (v1) at (0,-1.15);
  \coordinate (vj) at (0,-.6);
  \coordinate (u) at (0,1.5);
  \coordinate (x) at (-.5,.5);
  \coordinate (t) at (0,0);
  \coordinate (z) at (0,1);
  \coordinate (z1) at (.5,0.5);

  
  \draw[lf edge] (p) arc[start angle=180,end angle=90,radius=2]
      arc[start angle=90,end angle=0,radius=2]
      arc[start angle=0,end angle=-90,radius=2];
  \draw[lf fan] (p) arc[start angle=180,end angle=270,radius=2];
  \draw[lf fan] (s)--(v0)--(v1);
  \draw[lf fan] (p)--(v0) (p)--(v1) (p)--(vj);
  \draw[lf fan,dash pattern=on 0pt off 3pt,line cap=round] (v1)--(vj);
  \draw[draw=black,line width=.65pt]
      (0,.5) circle[radius=.5];
  \node at (0,.5) {$B$};
  \draw[lf edge] (p)--(x);
  \draw[lf edge] (vj)--(t);

  \draw[lf edge] (p)--(u);
  \ifnum#1=1
    \draw[lf edge] (r).. controls (0.5,1.75) and (0.5,1.25) .. (z) (u)--(z);
  \else
    \draw[lf edge] (r).. controls (0.5,1.75) and (0.5,1.25) .. (z1) (u)--(z);
  \fi
  \draw[lf edge,preaction={draw=white,line width=3.2pt}] (u)--(q);

  \foreach \a in {p,s,v0,v1,vj} \node[lf Hvertex] at (\a) {};
  \foreach \a in {q,r,u,x,t,z} \node[lf vertex] at (\a) {};
  \node[left=3pt] at (p) {$p$};
  \node[above=3pt] at (r) {$r$};
  \node[right=3pt] at (q) {$q$};
  \node[below=3pt] at (s) {$s$};
  \node[above left=0pt] at (u) {$u$};
  \node[below left=2pt] at (x) {$x$};
  \node[below left=1pt] at (t) {$e$};
  \node[right=3pt] at (v0) {$v_0$};
  \node[right=3pt] at (v1) {$v_1$};
  \node[right=3pt] at (vj) {$v_j$};
  \node[text=black] at (-1.15,-1.2) {$H$};

  \node at (1.7,1.8) {$P$};
  \node at (-1.7,1.8) {$C$};
  \ifnum#1=1
    \node[above left=0pt] at (z) {$z$};
    \node at (0,-2.85) {(a) $D_H=\{z\}$};
  \else
    \node[lf vertex] at (z1) {};
    \node[right=2pt] at (z1) {$z_1$};
    \node[above left=0pt] at (z) {$z_2$};
    \node at (0,-2.85) {(b) $D_H=\{z_1,z_2\}$};
  \fi
\end{tikzpicture}%
}
\centering
\begin{minipage}[t]{.48\linewidth}
\centering\lfPanel{1}
\end{minipage}\hfill
\begin{minipage}[t]{.48\linewidth}
\centering\lfPanel{2}
\end{minipage}
\endgroup

\caption{The configuration for a terminal fan of type~\textup{(i)} and $j=2$. The bold edges form $H$ and $P=C-V(H)=rq$. \textup{(a)} corresponds to $|D_H|=1$ and \textup{(b)} corresponds to $|D_H|\geq 2$.}
\label{fig:terminal-attachments}
\end{figure}
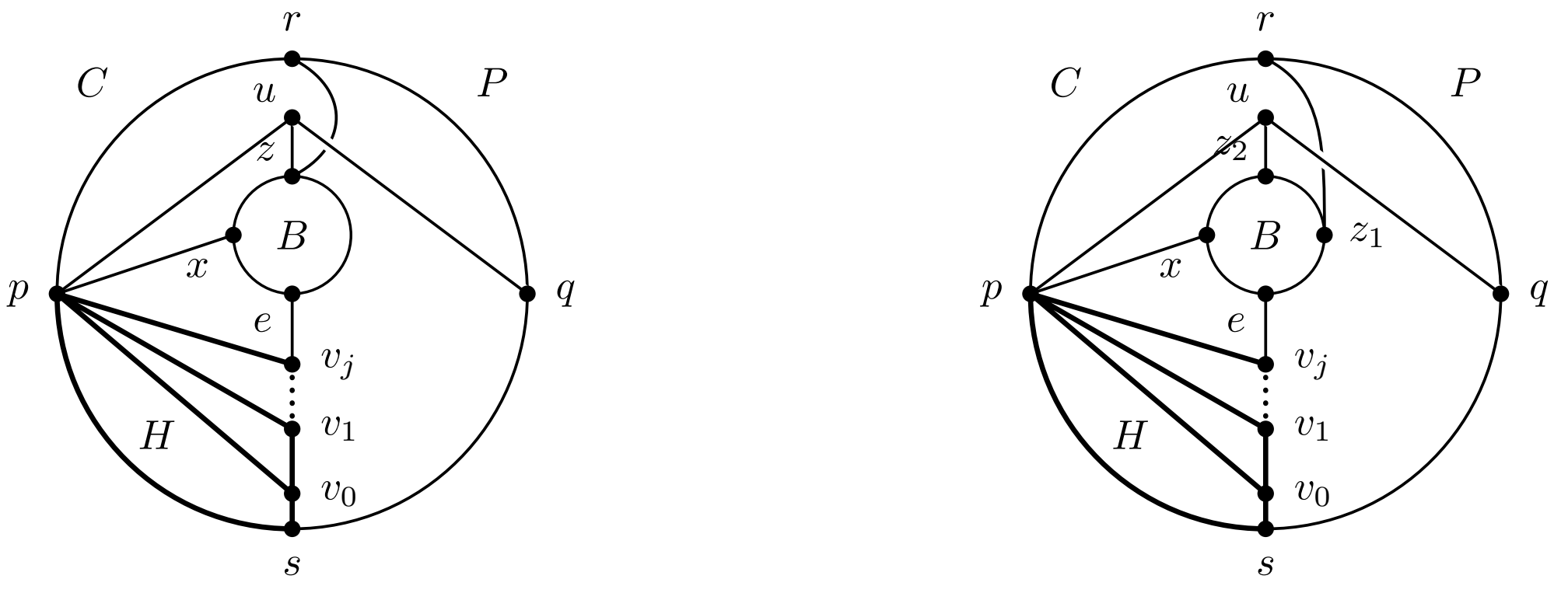
\begin{proof}
Suppose $H$ has centre $p$. Define $R:=T-V(H)$ and $P:=C-V(H)$. The graph $R$ is a tree containing $b$, and $P$ is a path on two or three vertices. If $D_H=\varnothing$, all neighbours of $B$ outside $B$ lie in $\{p,v_j\}$, contradicting 3-connectivity. Thus $D_H$ is nonempty.

\begin{enumerate}[label=\textup{(\arabic*)},leftmargin=*]
\item \emph{If $|D_H|=1$, then $H$ is not contractable.}

Suppose $D_H=\{z\}$, then deleting $z$ separates $B-z$ from $P$ in $G-H$ (see Figure 3 (a)), so $H$ is not contractable by Lemma \ref{lem:contract}.

\item \emph{If $|D_H|\ge2$, then $H$ is contractable.}

By Lemma \ref{lem:contract}, it suffices to show $G-H$ is 2-connected. By Lemma \ref{lem:block}(i), it suffices to show $(G-H)/B$ is 2-connected since $|D_H|\ge2$. By Lemma \ref{lem:treepath}, it suffices to show every end of $P$ meets $R$ and every leaf of $R$ meets $P$.

Every end of $P$ in $\{r,s\}$ meets $R$ since every vertex in $V(S\cap F)$ has its cycle neighbours in $V(H)\cup\{p,q\}$ and hence not adjacent to that end. If $q$ is an end of $P$, the fan is of type (i), so $u\notin S$ and $q$ meets $u\in R$ (see Figure 3 (b)). Thus every end of $P$ meets $R$.

Observe that every vertex in $\{r,s\}\cap V(H)$ has its unique $F$-neighbour in $S$. An $R$-old leaf of $R$ meets $P$: it has two cycle neighbours and cannot meet a vertex of $C\cap H$ except at $p$. An $R$-new leaf $v\ne b$ has degree 2 in $T$ and hence has at least one neighbour on $C$. By the maximality of $S$, $vp$ is not an edge and hence $v$ meets $P$.

It remains to consider $b$ when it is a leaf of $R$. If $R=\{b\}$, the endpoint argument already implies $b$ meets $P$. Thus we may assume $|R|\geq 2$. Since $b$ is a leaf, at most one vertex of $B$ meets $R-b$. Since $|D_H|\ge2$, some vertex of $B$ meets $P$ and hence $b$ meets $P$. It follows that $(G-H)/B$ is 2-connected as desired.

\item If $H$ is not contractable, then $v_j$ has a neighbour in $B$. 

By (2) and $D_H\neq \varnothing$, we have $|D_H|=1$ and write $D_H=\{z\}$. If $v_j$ has no neighbour in $B$, then $p$ is the only vertex of $H$ that may be adjacent to $B$. Thus $\{p,z\}$ separates $B-z$ from $P$ in $G$, a contradiction.
\end{enumerate}
\end{proof}

\begin{lemma}\label{lem:crown}
Let $e\in E(T)$, $S$ be the component of $T-e$ disjoint from $b$, $w$ be the end of $e$ in $S$ and $U:=V(S)\setminus\{w\}$. If every $v\in U$ is a leaf of $T$ with $N_G(v)=\{p,q,w\}$ and $|U|+|Z_w|\ge2$, then $G$ has a contractable fan, a contractable 2-crown, or a contractable 3-crown.

\end{lemma}

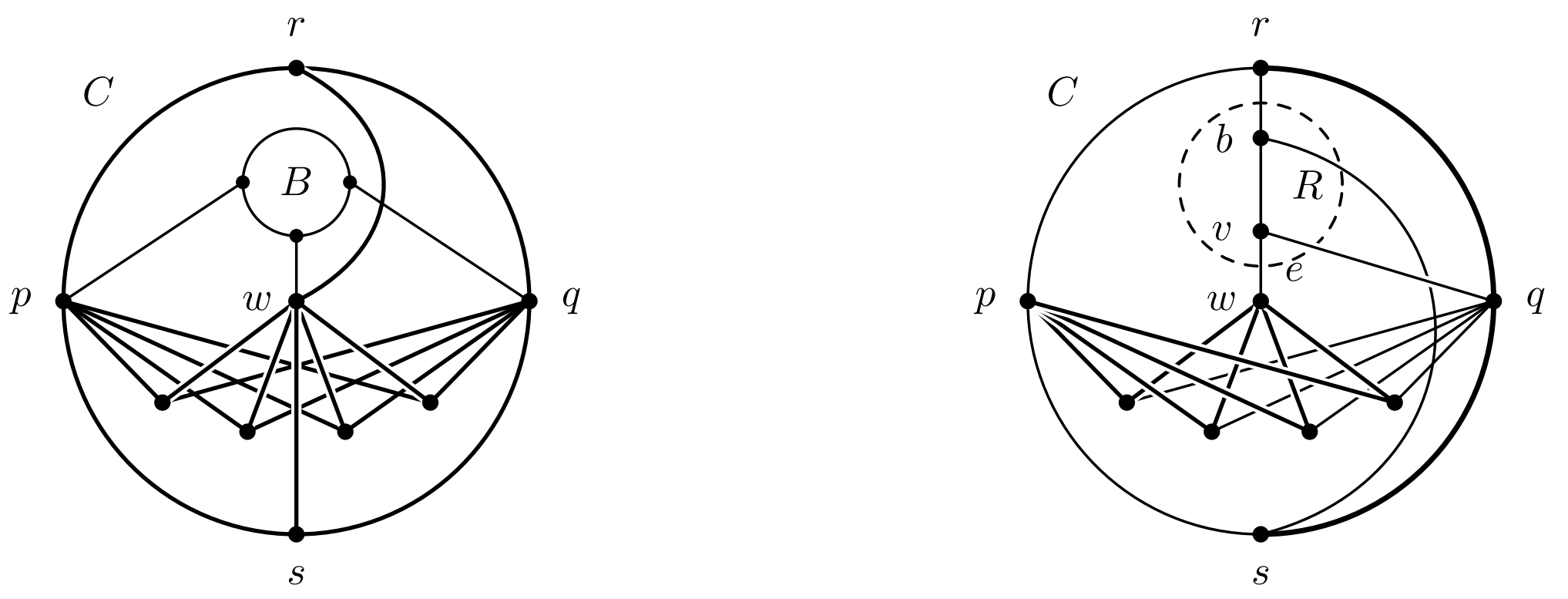
\begin{figure}[H]
\centering
\begin{minipage}[t]{.48\linewidth}
\centering
\begin{tikzpicture}[x=1cm,y=1cm,line cap=round,line join=round,
  every node/.style={font=\small},
  cv/.style={circle,fill=black,inner sep=0pt,minimum size=4pt},
  ch/.style={cv,fill=black},
  ce/.style={draw=black,line width=.65pt},
  he/.style={draw=black,line width=1pt},
  lab/.style={fill=white,inner sep=1.2pt}]
  \path[use as bounding box] (-2.65,-2.65) rectangle (2.65,2.55);
  \coordinate (p) at (-2,0); \coordinate (r) at (0,2);
  \coordinate (q) at (2,0); \coordinate (s) at (0,-2);
  \coordinate (w) at (0,0);
  \coordinate (u1) at (-1.15,-.87);
  \coordinate (u2) at (-.42,-1.12);
  \coordinate (u3) at (.42,-1.12);
  \coordinate (u4) at (1.15,-.87);
  \coordinate (xp) at (-.46,1.02);
  \coordinate (xq) at (.46,1.02);
  \coordinate (xw) at (0,.56);
  \draw[he] (0,0) circle[radius=2];
  \foreach \v in {u1,u2,u3,u4}{\draw[he] (p)--(\v)--(q);}
  \foreach \v in {u1,u2,u3,u4}{
    \draw[he,preaction={draw=white,line width=2.5pt}] (w)--(\v);}
  \draw[he,preaction={draw=white,line width=2.5pt}] (w)--(s);
  \draw[he,preaction={draw=white,line width=2.5pt}]
    (w) .. controls (1,0.5) and (1,1.5) .. (r);
  \draw[ce] (p)--(xp) (q)--(xq) (w)--(xw);
  \draw[draw=black,line width=.65pt,]
    (0,1.02) circle[radius=.46];
  \node at (0,1.03) {$B$};
  \foreach \v in {xp,xq,xw}{\node[cv,minimum size=3.4pt] at (\v) {};}
  \foreach \v in {p,r,q,s,w,u1,u2,u3,u4}{\node[ch] at (\v) {};}
  \node[left=4pt] at (p) {$p$}; \node[above=4pt] at (r) {$r$};
  \node[right=4pt] at (q) {$q$}; \node[below=4pt] at (s) {$s$};
  \node[left=2pt] at (w) {$w$};
  \node at (-1.7,1.8) {$C$};
\end{tikzpicture}
\par\smallskip
{\small (a) $R=\{b\}$, $Z=\{r,s\}$.}
\end{minipage}\hfill
\begin{minipage}[t]{.48\linewidth}
\centering
\begin{tikzpicture}[x=1cm,y=1cm,line cap=round,line join=round,
  every node/.style={font=\small},
  cv/.style={circle,fill=black,inner sep=0pt,minimum size=4pt},
  ch/.style={cv,fill=black},
  ce/.style={draw=black,line width=.65pt},
  he/.style={draw=black,line width=1pt},
  lab/.style={fill=white,inner sep=1.2pt}]
  \path[use as bounding box] (-2.65,-2.65) rectangle (2.65,2.55);
  \coordinate (p) at (-2,0); \coordinate (r) at (0,2);
  \coordinate (q) at (2,0); \coordinate (s) at (0,-2);
  \coordinate (w) at (0,0); \coordinate (v) at (0,.6);
  \coordinate (b) at (0,1.4);
  \coordinate (u1) at (-1.15,-.87);
  \coordinate (u2) at (-.42,-1.12);
  \coordinate (u3) at (.42,-1.12);
  \coordinate (u4) at (1.15,-.87);
  \draw[ce] (0,0) circle[radius=2];
  \draw[black,line width=1.3pt] (r) arc[start angle=90,end angle=-90,radius=2];
  \draw[draw=black,dashed,line width=.7pt,]
    (0,1) circle[radius=.7];
  \draw[ce] (r)--(b)--(v)--(w);
  \draw[ce] (s) .. controls (2,-1.5) and (2,1) .. (b);
  \draw[ce,preaction={draw=white,line width=2.4pt}] (q)--(v);
  \foreach \u in {u1,u2,u3,u4}{\draw[ce] (q)--(\u);}
  \foreach \u in {u1,u2,u3,u4}{
    \draw[he,preaction={draw=white,line width=2.5pt}] (p)--(\u)--(w);}
  \foreach \n in {r,q,s,v,b}{\node[cv] at (\n) {};}
  \foreach \n in {p,w,u1,u2,u3,u4}{\node[ch] at (\n) {};}
  \node[left=4pt] at (p) {$p$}; \node[above=4pt] at (r) {$r$};
  \node[right=4pt] at (q) {$q$}; \node[below=4pt] at (s) {$s$};
  \node[left=2pt] at (w) {$w$};
  \node[left=3pt] at (v) {$v$}; \node[left=3pt] at (b) {$b$};
  \node at (.4,1) {$R$};
  \node[right=2pt] at (0,.25) {$e$};
  \node at (-1.7,1.8) {$C$};
\end{tikzpicture}
\par\smallskip
{\small (b) $Z=\varnothing$, $H_p$ and $P_p=rqs$.}
\end{minipage}

\caption{The configurations in Lemma~\ref{lem:crown}. \textup{(a)} corresponds to \textup{(1)}. \textup{(b)} corresponds to \textup{(2)--(4)} in $G/B$ with $Z=\varnothing$. The bold edges in \textup{(b)} form $H_p$ and the complementary path $P_p$. The edge $e=wv$ separates $R$ from $S=\{w\}\cup U$.}
\label{fig:pendant-crowns}
\end{figure}
\begin{proof}
Let $Z:=Z_w$, and $R:=T-V(S)$. For any $k\in\{p,q\}$, let $H_k=G[\{k,w\}\cup U\cup Z]$ and $P_k=C-V(H_k)$. Every vertex of $U\cup Z$ is cubic with neighbourhood $\{p,q,w\}$. Since $|U|+|Z|\ge2$, each $H_k$ is a $2$-crown or a fan by Lemma~\ref{lem:forms}. Each $P_k$ is a nonempty path, and no vertex of $R$ meets $Z$. We first show that we can assume $|R|+|P_k|\ge3$ for any $k$.

\begin{enumerate}[label=\textup{(\arabic*)},leftmargin=*]
\item  \emph{Either $|R|+|P_k|\ge3$ for both choices of $k$ or $G$ has a contractable 3-crown.}

Suppose $|R|+|P_k|=2$ for any $k\in \{p,q\}$. Then $R=\{b\}$ and $Z=\{r,s\}$. Set $Y=U\cup\{r,s\}$ and $J=G-B$. Then $V(J)=\{p,q,w\}\cup Y$ and $N_G(y)=\{p,q,w\}$ for any $y\in Y$. See Figure 4 (a) for an example where $|U|=4$.

The specified leaf $u$ is in $U$: $R=\{b\}$ implies $u\in V(S)$ and $u\neq w$ since $w$ is adjacent to $r$ and $s$. Thus $|Y|\ge3$. Moreover, $pw$ is not an edge: otherwise both $p,w$ are noncubic on the triangle $pwrp$, contrary to Halin's lemma. Similarly we have $qw\notin E(G)$. Thus $J$ is an induced 3-crown. It is contractable since its complement is $B$. 

In view of (1), we may assume $|R|+|P_k|\ge3$ for any $k\in \{p,q\}$ in the rest of the proof.

\item \emph{If a leaf $v$ of $R$ misses $P_k$, then $v$ is $R$-new and $k$ is the only neighbour of $v$ on $C$.}

If $v$ is $R$-old, then it has two neighbours in $C$, none in $Z$, and therefore meets any $P_k$. Thus $v$ is $R$-new. No vertex of $R$ meets $Z$ and notice that $V(C-P_k)=\{k\}\cup Z$. Thus the only possible neighbour of $v$ on $C$ is $k$ since $v$ misses $P_k$ (see Figure 4 (b)). Notice that $d_T(v)=2$ and $v$ has at least one neighbour on $C$. It follows that $k$ is the only neighbour of $v$ on $C$.

\item \emph{Either $H_p$ or $H_q$ is contractable in $G/B$.}

Since $|R|+|P_k|\ge3$ and by Lemma~\ref{lem:treepath}, it suffices to show that there exists $k\in \{p,q\}$ such that every leaf of $R$ meets $P_k$ and every end of $P_k$ meets $R$. For any $k\in \{p,q\}$, every end of $P_k$ in $\{r,s\}$ is not in $Z$ and hence has its $F$-neighbour in $R$.

We first consider the case where every leaf of $R$ meets $P_p$ and $P_q$. Then $H_k$ is contractable if some $P_k$ has both ends in $\{r,s\}$ and we are done. If this is not the case, then either $r$ or $s$ is in $Z$ and hence not adjacent to a vertex in $R$. Note that $R$ contains at least one $R$-old leaf, which has at least one neighbour, say $q$, in $\{p,q\}$. Then both ends of $P_p$ meet $R$ and $H_p$ is contractable.

It remains to consider the case where some leaf $v$ misses some $P_k$. We may assume $k=q$. Then $v$ is an $R$-new leaf and $q$ is the only neighbour of $v$ on $C$ by (2) (also see Figure 4 (b)). It follows that every leaf of $R$ meets $P_p$, and $H_p$ is contractable if $P_p$ has both ends in $\{r,s\}$. If this is not the case, then $q$ is an end of $P_p$. Now both ends of $P_p$ still meet $R$ since $q$ is adjacent to $v$, and therefore $H_p$ is still contractable.

\item \emph{Either $H_p$ or $H_q$ is contractable in $G$.}

We may assume $H_p$ is contractable in $G/B$ by (3). There is an edge between $B$ and $(R-b)\cup(\{r,s\}\setminus Z)$: The statement is clear if $R-b$ is nonempty, and if $R-b$ is empty, then $\{r,s\}\setminus Z$ is nonempty by (1), and a vertex in it meets $B$. Thus $B$ has an attachment in both $G-H_p$ and $G-H_q$.

If $H_p$ is not contractable in $G$, the attachment set of $B$ in $G-H_p$ is a singleton $\{z\}$ by Lemma~\ref{lem:block}(i). Then $w$ has a neighbour $x$ in $B$ since otherwise $\{p,z\}$ would be a 2-vertex-cut of $G$. In particular, $e=wb$. It follows that $p$ has a neighbour $y\notin \{x,z\}$ in $B$ since otherwise $\{x,z\}$ would be a vertex-cut of size at most 2.

Since $e=wb$, the only possible $R$-new leaf is $b$. A similar leaf and endpoint argument in (3) implies that $H_q$ is contractable in $G/B$. Moreover, the common attachment in both $G-H_p$ and $G-H_q$ is $z$ since the attachment set of $B$ in $G-H_p$ is a singleton. On the other hand, $y$ is also an attachment in $G-H_q$. Therefore $G-H_q$ is 2-connected by Lemma~\ref{lem:block}(i).
\end{enumerate}
\end{proof}

\begin{lemma}\label{lem:quad}
In the quadrilateral setup above, $G$ has a contractable fan, a contractable $2$-crown, or a contractable 3-crown.
\end{lemma}
\begin{proof}
Suppose none exists. Root $T$ at $b$, and let $L$ consist of the vertices $v$ such that $\{v\}$ is a component of $T-b$ and $N_C(v)=\{p,q\}$. We say a component $H$ of $T-b$ gives a terminal fan if $G[V(H)\cup V(C)]$ contains a terminal fan.

\begin{enumerate}[label=\textup{(\arabic*)},leftmargin=*]
\item \emph{Every component of $T-b$ either gives a terminal fan containing the whole component, or contains a unique leaf in $L$.}

By Lemma~\ref{lem:fan}, if a component of $T-b$ gives a terminal fan, its final vertex meets $B$. Therefore the terminal fan contains the whole component of $T-b$. Let $Q$ be a component of $T-b$ giving no terminal fan and $v$ be a leaf of $T$ in $Q$. 

If $v$ has cycle neighbours in $\{r,s\}$, then $G$ has a contractable configuration by setting $S=\{v\}$ in Lemma~\ref{lem:crown}, a contradiction. If $v$ has cycle neighbours in both $\{p,q\}$ and $\{r,s\}$, then $Q$ gives a terminal fan, a contradiction. Thus $v$ has cycle neighbours in $\{p,q\}$. 

If $|Q|\ge2$, choose $w\in Q$ farthest from $b$ among vertices having a child. Its children form a nonempty set $U$ of leaves, all with neighbourhood $\{p,q,w\}$, and $S=\{w\}\cup U$ is cut off by the edge towards $b$. Lemma~\ref{lem:crown} forces $|U|+|Z_w|<2$. Hence $U=\{v\}$ and $Z_w=\varnothing$. Now $d_T(w)=2$, and $Q$ contains a path defining a terminal fan of type (ii), a contradiction.

\item \emph{At most one component gives a terminal fan.}

Suppose two components give terminal fans $H_1,H_2$ and let $D_i=N_B(G-V(B)-V(H_i))$. By Lemma~\ref{lem:fan}, both $D_i$ are singletons. If $H_1,H_2$ have different centres, then they are disjoint: their intersection is contained in $\{r,s\}$, while both vertices are cubic and cannot have neighbours in both defining components. Every attachment vertex of $B$ then belongs to $D_1\cup D_2$, contrary to the attachment inequality~\eqref{eq:attachments}.

Thus $H_1,H_2$ have the same centre, say $p$. We claim that $D_1\cap D_2\neq\varnothing$. By definition of terminal fan, each $H_i$ contains at most one of $u,r,s$. Thus some $v\in\{u,r,s\}$ lies outside both $H_i$'s. If $v=u$, its component gives an attachment in $D_1\cap D_2$. If $v\in\{r,s\}$, its unique $F$-neighbour cannot lie on either fan path since a vertex in any $F\cap H_i$ has its cycle neighbours contained in $V(H_i)\cup\{p,q\}$. That neighbour therefore lies in $B$ or in another component, again giving an attachment in $D_1\cap D_2$, proving the claim. 

Each $D_i$ is a singleton, so we may assume $D_1=D_2=\{z\}$. Every edge from $B-z$ to the outside of $B$ now ends at $p$, so $\{p,z\}$ is a separator, a contradiction.

\item \emph{The remaining vertices induce a contractable fan or $2$-crown.}
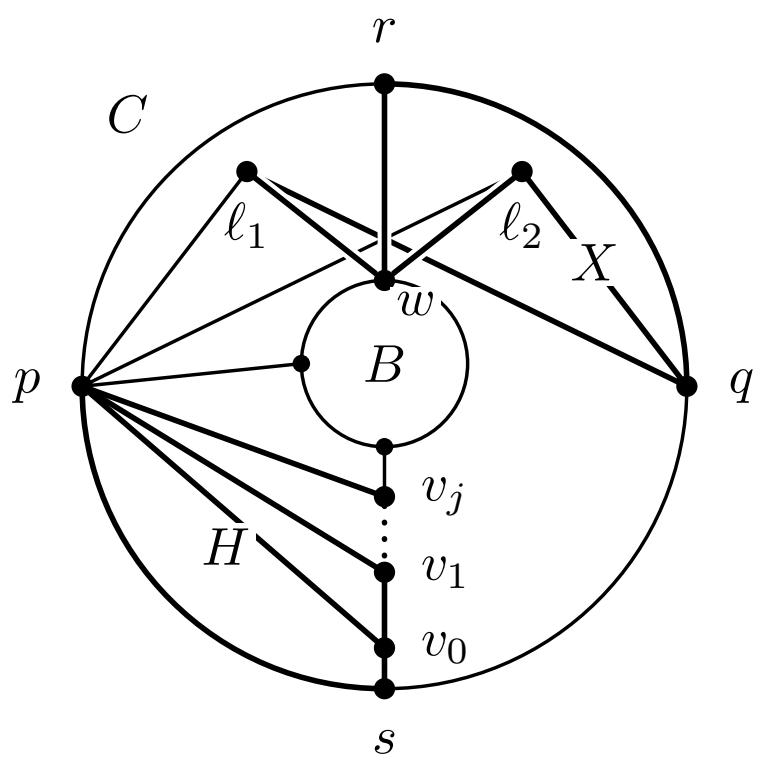
\begin{figure}[H]
  \centering
  \begin{tikzpicture}[x=1cm,y=1cm,line cap=round,line join=round,
    every node/.style={font=\small},
    cv/.style={circle,fill=black,inner sep=0pt,minimum size=4pt},
    hv/.style={cv,fill=black},
    xv/.style={cv,fill=black},
    ce/.style={draw=black,line width=.65pt},
    he/.style={draw=black,line width=1.05pt},
    xe/.style={draw=black,line width=1.05pt},
    lab/.style={fill=white,inner sep=1pt}]
    \path[use as bounding box] (-2.8,-2.6) rectangle (2.8,2.6);
    \coordinate (p) at (-2,0); \coordinate (r) at (0,2);
    \coordinate (q) at (2,0); \coordinate (s) at (0,-2);
    \coordinate (w) at (0,.70);
    \coordinate (l1) at (-.91,1.42); \coordinate (l2) at (.91,1.42);
    \coordinate (vj) at (0,-.73); \coordinate (v1) at (0,-1.23);
    \coordinate (v0) at (0,-1.73);
    \coordinate (t) at (0,-.4); \coordinate (x) at (-.55,.15);
    
    \draw[ce] (0,0) circle[radius=2];
    \draw[he] (p) arc[start angle=180,end angle=270,radius=2];
    \draw[xe] (r) arc[start angle=90,end angle=0,radius=2];
    \draw[ce] (p)--(l1) (p)--(l2) (p)--(x);
    \draw[xe] (q)--(l1) (q)--(l2);
    \draw[xe,preaction={draw=white,line width=2.6pt}] (l1)--(w)--(l2);
    \draw[xe,preaction={draw=white,line width=2.6pt}] (r)--(w);
    \draw[draw=black,line width=.65pt,]
      (0,.15) circle[radius=.55];
    \node at (0,.15) {$B$};
    \draw[ce] (t)--(vj);
    \draw[he] (v1)--(v0)--(s);
    \draw[he,dash pattern=on 0pt off 3pt,line cap=round] (v1)--(vj);
    \foreach \v in {vj,v1,v0}{\draw[he] (p)--(\v);}
    \foreach \v in {p,s,vj,v1,v0}{\node[hv] at (\v) {};}
    \foreach \v in {q,r,w,l1,l2}{\node[xv] at (\v) {};}
    \foreach \v in {t,x}{\node[cv,minimum size=3.3pt] at (\v) {};}
    \node[left=4pt] at (p) {$p$}; \node[above=4pt] at (r) {$r$};
    \node[right=4pt] at (q) {$q$}; \node[below=4pt] at (s) {$s$};
    \node[below right=1pt,lab] at (w) {$w$};
    \node[below=2pt] at (l1) {$\ell_1$};
    \node[below=2pt] at (l2) {$\ell_2$};
    \node[right=3pt] at (vj) {$v_j$};
    \node[right=3pt] at (v1) {$v_1$};
    \node[right=3pt] at (v0) {$v_0$};
    \node[black,lab] at (-1.05,-1.06) {$H$};
    \node[black,lab] at (1.39,.82) {$X$};
    \node at (-1.7,1.8) {$C$};
  \end{tikzpicture}
  
  \caption{The configuration in Lemma \ref{lem:quad} (3). The bold edges form the fan $H$ and the $2$-crown $X$, where $Y=\{r,\ell_1,\ell_2\}$. The vertex $w$ is the unique member of $D_H$.}
  \label{fig:final-crown}
  \end{figure}

If there is no terminal fan, then $G-B$ is a contractable $2$-crown. Thus we suppose there is a unique terminal fan $H$ with centre $p$. Let $D_H=\{w\}$ and $Y=(L\cup\{r,s\})\setminus V(H)$. See Figure 5 for an example where $|Y|=3$. Then each vertex in $Y$ has its unique $F$-neighbour $w\in V(B)$, and hence every vertex of $Y$ is cubic with neighbourhood $\{p,q,w\}$. Note that at least two of $u,r,s$ are outside $H$ and if $u$ is outside, then it belongs to $L$ by (1) and (2). Thus $|Y|\ge2$ and $X:=G[\{q,w\}\cup Y]$ is a $2$-crown or a fan by Lemma \ref{lem:forms}.

Every edge from $B-w$ to $G-B$ ends in $H$, and $V(G-X)=V(H)\cup(V(B)\setminus\{w\})$.
Thus $G-X$ is 2-connected by Lemma \ref{lem:block} (ii), and $X$ is contractable by Lemma \ref{lem:contract}, a contradiction.
\end{enumerate}
\end{proof}

\begin{proof}[Proof of Theorem \ref{thm:main}]
Let $G$ be minimally 3-connected and neither a wheel nor a graph $K_{3,m}$. By Theorem \ref{thm:dirac}, $G$ contains two vertex-disjoint cycles. Choose $C,B$ by Theorem \ref{lem:core}, with $B$ 2-connected as explained after its statement. Alternative \textup{(a)} gives a contractable induced cycle, and alternative \textup{(d)} contradicts Halin's lemma. Alternative \textup{(c)} is handled by Lemma \ref{lem:triangle},
and alternative \textup{(b)} by Lemma \ref{lem:quad}.
\end{proof}

\section{Concluding remarks}\label{sec:terminal}

It is not difficult to show that if no 2-crown contained in the 3-crown is contractable, then there are exactly 3 edges between that 3-crown and $B$, forming a matching 3-edge cut, which is a 3-edge cut forming a matching. In fact, let $J$ be a contractable 3-crown in a minimally 3-connected graph $G$ with centres $x_1,x_2,x_3$. Since $J$ is contractable, the 2-crown $J-x_i$ is contractable exactly when $x_i$ has at least two neighbours in $G-J$. If none is contractable, then each $x_i$ has exactly one neighbour $x_i'$ in $G-J$. Then $x_1'$, $x_2'$, $x_3'$ are distinct since $G$ is 3-connected. Thus $x_1x_1'$, $x_2x_2'$, $x_3x_3'$ form a matching 3-edge cut.  

Note that a matching 3-edge cut naturally gives a decomposition of a 3-connected graph into smaller 3-connected graphs. Therefore, our result can be restated in the language of $Y$-sum: if $G_1,G_2$ have cubic vertices $v_1,v_2$, a \emph{$Y$-sum} of $G_1$ and $G_2$ is obtained from the disjoint union of $G_1-v_1$ and $G_2-v_2$ by joining the three neighbours of $v_1$ to the three neighbours of $v_2$ by a matching of size three. Note that a $Y$-sum of two 3-connected graphs is again 3-connected. A 3-connected graph $G$ is said to be $Y$-$reducible$ if it is a $Y$-sum of two smaller 3-connected graphs. Then a 3-connected graph having a matching 3-edge cut is $Y$-reducible and the following corollary is immediate.

\begin{corollary}\label{coro}
  Let $G$ be a minimally 3-connected graph. At least one of the following holds.
  \begin{enumerate}[label=\textup{(\arabic*)}]
  \item The graph $G$ is $Y$-reducible.
  \item The graph $G$ has a contractable induced cycle, a contractable fan,
        or a contractable 2-crown.
  \item $G\in\Kcal\cup\Wcal$.
\end{enumerate}
\end{corollary}

Using the above corollary, our result can be conversely translated into a construction theorem for 3-connected graphs: each 3-connected graph can be obtained from a wheel or $K_{3,n}$ by repeatedly adding an edge between nonadjacent vertices, taking $Y$-sum, and replacing a vertex $v$ by a configuration $H$ that is either an induced cycle, an induced fan or an induced 2-crown and add joining edges such that:

(1) Any vertex originally adjacent to $v$ is joined to at least one vertex in $H$.

(2) Every vertex in $H$ has degree at least 3 and for every $v\in V(H)$, $H-v$ has at least two distinct neighbours outside $H$.

It is not difficult to see that they are necessary and sufficient for the expanded graph to be 3-connected. There is no direct anology for minimally 3-connectedness since the vertex expansion operation may destroy minimally 3-connectedness. This construction differs from the classical chain theorems in that the size increment in each step is not bounded but the increment structure is guaranteed to be 2-connected.

\end{document}